\documentclass[12pt]{amsart}
\usepackage{amscd,amssymb}
\usepackage[graph,frame,poly,arc]{xy}  
\usepackage[plainpages,backref,urlcolor=blue]{hyperref}

\theoremstyle{plain}
\newtheorem{thm}[subsection]{Theorem}
\newtheorem{lem}[subsection]{Lemma}
\newtheorem{prop}[subsection]{Proposition}
\newtheorem{cor}[subsection]{Corollary}

\theoremstyle{definition}
\newtheorem{rk}[subsection]{Remark}
\newtheorem{definition}[subsection]{Definition}

\newtheorem{conj}[subsection]{Conjecture}

\numberwithin{equation}{section}
\newcommand{\A}{{\mathcal A}}
\newcommand{\B}{{\mathcal B}}

\newcommand{\al}{{\alpha}}

\newcommand{\Z}{\mathbb{Z}}

\newcommand{\C}{\mathbb{C}}
\newcommand{\K}{\mathbb{K}}
\newcommand{\PP}{\mathbb{P}}

\DeclareMathOperator{\mult}{mult}

\DeclareMathOperator{\Der}{Der}

\begin{document}

\title [On supersolvable line arrangements]
{On complex supersolvable line arrangements}

\author[Takuro Abe]{Takuro Abe$^1$}
\address{Institute of Mathematics for Industry,
Kyushu University,
Fukuoka 819-0395, Japan.}
\email{abe@imi.kyushu-u.ac.jp}

\author[Alexandru Dimca]{Alexandru Dimca$^2$}
\address{Universit\'e C\^ ote d'Azur, CNRS, LJAD, France and Simion Stoilow Institute of Mathematics,
P.O. Box 1-764, RO-014700 Bucharest, Romania}
\email{dimca@unice.fr}


\thanks{
$^1$This work is 
partially supported by KAKENHI, JSPS Grant-in-Aid for Scientific Research (B) 16H03924.\\
\ \ $^2$ This work has been supported by the French government, through the $\rm UCA^{\rm JEDI}$ Investments in the Future project managed by the National Research Agency (ANR) with the reference number ANR-15-IDEX-01 and by the Romanian Ministry of Research and Innovation, CNCS - UEFISCDI, grant PN-III-P4-ID-PCE-2016-0030, within PNCDI III}

\subjclass[2010]{Primary 14H50; Secondary  14B05, 13D02, 32S22}

\keywords{supersolvable line arrangement, double points }

\begin{abstract} We show that the number of lines in an 
$m$--homogeneous supersolvable line arrangement is upper bounded by $3m-3$ and we classify the $m$--homogeneous supersolvable line arrangements with two modular points up-to lattice-isotopy.
A lower bound for the number of double points $n_2$ in an $m$--homogeneous supersolvable line arrangement of $d$ lines is also considered.
 When $3 \leq m \leq 5$, or when $m \geq \frac{d}{2}$, or when there are at least two modular points, we show that $n_2 \geq \frac{d}{2}$, as conjectured by B. Anzis and S. O. Toh\u aneanu. This conjecture is shown to hold also for supersolvable line arrangements obtained as cones over generic line arrangements, or cones over  arbitrary line arrangements having a generic vertex.
\end{abstract}
 
\maketitle

\section{Introduction} 
Let $\A:f=0$ be a line arrangement  in the complex projective plane $\PP^2$.
An intersection point $p$ of $\A$ is called a {\it modular} point if for any other intersection point $q$ of $\A$, the line $\overline {pq}$ determined by the points $p$ and $q$ belongs to the arrangement $\A$. The arrangement $\A$ is {\it supersolvable} if it has a modular intersection point.
Supersolvable arrangements have many interesting properties, in particular they are free arrangements, see \cite{DHA,T1,T2} or  \cite[Prop 5.114]{OT} and \cite[Theorem 4.2]{JT}. For basic facts on free arrangements, please refer to \cite{DHA,OT, Yo}.
Recently, there is a renewed interest in supersolvable arrangements, see 
\cite{T1, CHMN, DMO,DStNSS, HaHa,T2}. 

It was shown by Hanumanthu and Harbourne in \cite[Corollary 2]{HaHa}, that supersolvable line arrangements having two modular points, say $p$ and $q$, with distinct multiplicity $m_p(\A) \ne m_q(\A)$ are easy to describe, being the union of two pencils with vertexes $p$ and $q$, with a common line $\overline {pq}$.
Hence it remains to study the {\it $m$-homogeneous supersolvable line arrangements, i.e. the supersolvable line arrangements in which all modular points have the same multiplicity} $m$. The case $m=2$ is easy, see \cite[Section 3.1]{HaHa}, and the $3$-homogeneous supersolvable line arrangements can be easily described, see \cite[Section 3.2.2]{HaHa} and Remark \ref{rkm=3},
so the really interesting cases are for $m \geq 4$.
The first result of our paper is the following.
\begin{thm}
\label{thm1}
Let $\A$ be a supersolvable arrangement in $\PP^2$. Let $p$ be a modular point of $\A$ with maximal multiplicity $m\geq 2$. Then $d:=|\A| \le 3m-3$. 
\end{thm}
In fact, this results holds for arrangements over any field $\K$ of characteristic zero, see Section 2 below.
For instance, when $m=100$, we have $d \leq 297$, and this gives a negative answer to the question about the existence of a 
$100$-homogeneous supersolvable line arrangement with $d\geq 7413$ in \cite[Remark 19]{HaHa}.
 The full monomial line arrangement 
 \begin{equation}
\label{FMA}
 \A(m-2,1,3): xyz(x^{m-2}-y^{m-2})(y^{m-2}-z^{m-2})(x^{m-2}-z^{m-2})=0,
\end{equation}  
for $m\geq 3$, has $d=3m-3$ and it is an $m$-homogeneous supersolvable line arrangement, and hence the bound in Theorem \ref{thm1} is sharp for this arrangement. $\A(m-2,1,3)$ has three modular points, located at $(1:0:0)$, $(0:1:0)$ and $(0:0:1)$ for $m>3$, while for $m=3$, there is a fourth modular point located at $(1:1:1)$.
The following result was obtained in \cite[Section 3]{HaHa}. We reprove this result as part of our Theorem \ref{thm1B} below, see also Remark \ref{rkHH} for the real supersolvable line arrangements.

\begin{thm}
\label{thmHH}
Let $\A$ be an $m$-homogeneous supersolvable line arrangement,
with $m\geq 3$. Then the number $M$ of modular points in $\A$
satisfies $M \leq 4$. Moreover,  any $m$-homogeneous supersolvable line arrangement $\A$ with $m\geq 3$, and having $M\geq 3$ modular points, is projectively equivalent to the arrangement $\A(m-2,1,3)$.
\end{thm}

To state our next result, we define first some subarrangements $\A$ in $\A(m-2,1,3)$ and consider the intersection lattices $L(\tilde \A)$ of the associated central plane arrangements $\tilde \A$ in $\C^3$. Set $n=m-2$ and denote by $\mu_n$ the multiplicative group of the $n$-th roots of unity. For $1 \leq k \leq n$, let 
$$W(n,k)=\mu_n^k \setminus \Delta,$$
where 
$$\Delta=\{(w_1,\ldots,w_k)  \ : \ w_i \in \mu_n \text{ and } w_{j_1}=w_{j_2} \text{ for some } j_1 \ne j_2\}.$$
For $w=(w_1,\ldots,w_k) \in W(n,k)$ we define the line arrangement
$$\A(w): f(w)=xyz(x^n-y^n)(x^n-z^n)\prod_{j=1,k}(z-w_jy)=0$$
in $\PP^2$ and the corresponding central plane arrangement $\widetilde {A(w)}:f(w)=0$ in $\C^3$. 
Let $\Sigma_k$ denote the symmetric group on $k$-elements, and consider the obvious action of  $\Sigma_k$ on the set $W(n,k)$. It is clear that
$f(w)=f(\sigma w)$, for any $\sigma \in \Sigma _k$.
The group $\mu_n$ also acts on the set $W(n,k)$ by translation, namely
$$aw=a(w_1,\ldots,w_k)=(aw_1,\ldots,aw_k)$$
for any $a \in \mu_n$. The coordinate change $y \mapsto ay$ shows that the line arrangements $\A(w)$ and $\A(aw)$ are projectively equivalent.
These two actions commute, and hence the direct product $N=\Sigma_k \times \mu_n$ acts on $W(n,k)$. Finally, the group $H=\mu_2$ acts on on the set $W(n,k)$ by
$$(-1)w=(-1)(w_1,\ldots,w_k)=(w_1^{-1},\ldots,w_k^{-1}).$$
The coordinate change $y \mapsto z$, $z \mapsto y$ shows that
the line arrangements $\A(w)$ and $\A((-1)w)$ are projectively equivalent. It follows that the semidirect product 
$$G=N \rtimes_{\phi}H$$
of $N$ and $H$, associated to the homomorphism $\phi: H \to Aut(N)$ given by 
$$\phi(-1)(\sigma,a)=(\sigma,a^{-1})$$
acts on the set $W(n,k)$. Moreover, to each class
$$[w] \in \overline{ W(n,k)}:=W(n,k)/G,$$
there is a line arrangement $\A(w)$ in $\PP^2$, well defined up-to projective equivalence, and a plane arrangement $\widetilde {A(w)}$ in $\C^3$, well defined up-to linear equivalence. For $k=0$ we define
$W(n,0)=\overline{ W(n,0)}$ to be a singleton $\{w_0\}$, and set
$$\A(w_0): f(w_0)=xyz(x^n-y^n)(x^n-z^n)=0.$$
Note that $\overline{W(n,n)}$ is also a singleton, and the corresponding
 line arrangement $\A(w)$  is precisely $\A(n,1,3)$. 
 
 For any $[w] \in \overline{ W(n,k)}$,
let $L(w)$ denote the intersection lattice of the plane arrangement $\widetilde {A(w)}$. Since the general linear group $GL_3(\C)$ is connected, it follows that 
for $[w]=[w']$, the plane arrangements $\widetilde {A(w)}$ and $\widetilde {A(w)}$
are lattice-isotopic, see \cite{RaLI1,RaLI2} or \cite[Definition 4.1]{DHA} for the notion of 
lattice-isotopic arrangements. This implies in particular that the two intersection lattices $L(w)$ and $L(w')$ are isomorphic, i.e. the arrangements $ \A(w)$ and $ \A(w')$ have the same combinatorics.
This class of line arrangements enters into the following result.
\begin{thm}
\label{thm1B}
Let $\A$ be an $m$-homogeneous supersolvable line arrangement with $m\geq3$, having $M\geq 2$ modular points. Then there is a unique integer $k$ such that $0 \leq k \leq m-2$ and a unique equivalence class
$[w] \in  \overline{ W(m-2,k)}$, such that the associated central arrangements $\tilde \A$ and $\widetilde {A(w)}$ are lattice-isotopic.
In particular, the two intersection lattices $L(\tilde \A)$ and $L(w)=L(\widetilde {A(w)})$ are isomorphic.
The number $M$ of modular points of the line arrangement $\A$ is given by
\begin{enumerate}

\item $M=2$ if $k<m-2$,

\item $M=3$ if $k=m-2 \geq 2$, and

\item $M=4$ if $k=m-2=1$.
\end{enumerate}
Moreover, when $k=m-2$, $\A$ is projectively equivalent to $ \A(m-2,1,3)$.
\end{thm}
Note that $|\A(m-2,1,3)|=3m-3$ and $|\A(w)|=2m-1+k <3m-3$ if $[w] \in  \overline{ W(m-2,k)}$.
Using Randell's results in \cite{RaLI1,RaLI2}, Theorem \ref{thmHH} and Theorem \ref{thm1B} imply the following. For basic information on the Milnor fiber $F(\A)$ and the monodromy transformation $h: F(\A) \to F(\A)$ we refer to \cite[Chapter 5]{DHA}.
\begin{cor}
\label{corTop}
Let $\A: f=0$ be an $m$-homogeneous supersolvable line arrangement with $m\geq3$, having at least two modular points. Then the combinatorics of $\A$ determines the diffeomorphism type of the complement $M(\A)=\PP^2 \setminus \cup_{L \in \A}L$ of $\A$ in $\PP^2$ and of the Milnor fiber $F(\A): f(x,y,z)=1$ in $\C^3$, as well as the monodromy operators $h^*:H^*(F(\A),\C) \to H^*(F(\A),\C)$.
\end{cor}

In order to study the remaining case of $m$-homogeneous supersolvable line arrangements, having one modular point, one notices first that such an arrangement is the cone $C(\A')_e$ over a simpler 
line arrangement $\A'$, constructed as follows. Assume that $\A'$ consists of $d'$ lines in $\PP^2$. Then we choose a point $p \in \PP^2$, the vertex of the cone, not situated on a line in $\A'$. We get $C(\A',p)_0$ by adding first to the lines in $\A'$, all the lines joining $p$ to the intersection points of $\A'$. Then we add $e \geq 0$ new lines through the vertex $p$, which meet $\A'$ only at simple points, and get in this way the line arrangement $C(\A',p)_e$.
To simplify the notation, we set $C(\A',p)_e=C(\A')_e$, even though the combinatorics of  $C(\A',p)_e$ depends heavily on $p$. We also  set $C(\A')=C(\A')_0$ and 
let $N'$ be the numbers of lines in $C(\A') \setminus \A'$.
In the special case when $\A'$ is a generic arrangement of $d'\geq 3$ lines, having 
$N={d' \choose 2}$ double points, such cones $C(\A')_e$  play a special role in the theory, see Remark \ref{rkEQ}.
If, in addition, $e=0$ and the point $p$ is chosen generic, we have $N'=N$.
In this  case one has for the arrangement $\A=C(\A')$: $d=d'+N$, $m=N\geq 3$, $n_2=N(d'-2)$,
$n_3=N$ and $n_m=1$, the other $n_j$ being zero. 
Note that for this arrangement $\A$ one has  $n_2 \geq d/2$.
Maybe this was the origin of the following conjecture, which occurs in \cite{T1}, where it was checked for all supersolvable line arrangements $\A$, not a pencil, which are either real, see \cite[Theorem 2.4]{T1}, or complex, satisfying $|\A| \leq 12$, see \cite[Corollary 3.5]{T1}.
\begin{conj}
\label{conj1}
Let $\A$ be a supersolvable line arrangement, consisting of $d=|\A|$ lines, which is not a pencil. Then $$n_2 \geq \frac{d}{2}.$$
\end{conj}
Theorem \ref{thmHH} and Theorem \ref{thm1B} imply the following result, since by deleting a line $z-ay=0$ with $a \in \mu_n$ from $\A(m-2,1,3)$, we lose one double point on the line $x=0$, but add at least $m-2$ double points on the deleted line. 
\begin{cor}
\label{cor1C}
Let $\A$ be an $m$-homogeneous supersolvable line arrangement, having at least two modular points. Then Conjecture \ref{conj1} holds for the line arrangement $\A$. 
\end{cor}
We prove, as a first step in the direction of understanding the $m$-homogeneous supersolvable line arrangements, having one modular point, the following three results.

\begin{thm}
\label{thm2B}
Let $\A$ be a free line arrangement, consisting of $d=|\A|$ lines and such that $\A $ is not a pencil. Let $m$ be the maximal multiplicity of an intersection point in $\A$.
If $2m\geq d$, then
$$n_2 \geq -2m^2+(3d-1)m-d^2+d \geq \frac{d}{2}.$$
In particular, Conjecture \ref{conj1} holds for an $m$-homogeneous supersolvable line arrangement $\A$, not a pencil, consisting of $d=|\A|$ lines such that $2m\geq d$.

\end{thm}

\begin{thm}
\label{thm2C}
Let $\A$ be a cone $C(\A')_e$ over a generic  arrangement $\A'$ of $d' \geq 3$ lines, with an arbitrary vertex $p$. Then 
Conjecture \ref{conj1} holds  for the supersolvable line arrangement $\A$.
\end{thm}
We can also treat the case of cones over arbitrary line arrangements $\A'$, but with a generic vertex.
\begin{prop}
\label{prop2C}
Let $\A$ be a cone $C(\A')_e$ over an arbitrary arrangement $\A'$ of $d' \geq 3$ lines, which is not a pencil, with a generic vertex $p$. Then 
 Conjecture \ref{conj1} holds for the supersolvable line arrangement $\A$.
\end{prop}

In \cite{HaHa}, the authors also put forth the following weaker conjecture.
\begin{conj}
\label{conj2}
Let $\A$ be a supersolvable line arrangement,  which is not a pencil. Then $$n_2 >0.$$
\end{conj}
Conjecture \ref{conj2} is proved for the 
$m$--homogeneous supersolvable line arrangements, with $3 \leq m \leq 4$, see \cite[Theorem 17]{HaHa}. Theorem \ref{thm1} implies the following stronger result.

\begin{cor}
\label{corthm2}
Any $m$--homogeneous supersolvable line arrangement which is not a pencil satisfies Conjecture \ref{conj1} if $3 \leq m \leq 5$.
\end{cor}

The proofs of the main results use a variety of techniques. For instance, the multiarrangements as in \cite{A1,A2,Z} are used in the proof of Theorem \ref{thm1}.
Fine algebraic properties of the Jacobian module of a generic line arrangement, and of its saturation, due to Rose and Terao \cite{RT}, Yuzvinsky \cite{Yu}, and Hassanzadeh and Simis \cite{HS}, are used in the proof of Theorem \ref{thm2C}. It may be a challenging question to find a purely combinatorial proof for this latter result.

\section{Proof of Theorem \ref{thm1}} 

We start by recalling some general preliminaries needed in the proof. For a general reference, see \cite{OT} and \cite{Yo}. In this section $\K$ is a field of characteristic zero. Note that an arrangement $\A$ in $\PP^2$ is the same as a central hyperplane arrangement in $\K^3$. Let 
$V=\K^3,\ S:=\mbox{Sym}^*(V^*)=\K[x,y,z]$ and 
$\Der S:=S \partial_{x} \oplus S \partial_{y}\oplus S \partial_{z}$. 
For $H \in \A$, we fix one 
linear form $\alpha_H \in V^*$ with $\ker \alpha_H=H$. The \textbf{logarithmic 
derivation module} $D(\A)$ is defined by 
$$
D(\A):=\{
\theta \in \Der S \mid \theta(\alpha_H) \in S\alpha_H\ (\forall H \in \A)\}.
$$
In general $D(\A)$ is a reflexive, but not a free $S$-module. We say that $\A$ is 
\textbf{free} with $\exp(\A)=(1,d_2,d_3)$ if 
\begin{equation}
\label{FA1}
D(\A) \simeq S[-1] \oplus S[-d_2] \oplus S[-d_3]
\end{equation}
as a graded $S$-module. If $\A$ is supersolvable with a modular point 
$p$ of multiplicity $m$, then it is known that 
$\A$ is free with 
\begin{equation}
\label{FA2}
\exp(\A)=(1,m-1,|\A|-m),
\end{equation} see Theorem 4.58 in \cite{OT} for 
example. Next let us introduce the Ziegler restriction that enables us to 
construct a multiarrangement from a simple arrangement canonically. 

\begin{definition}[\cite{Z}]
Let $H \in \A$. Then the \textbf{Ziegler restriction} $(\A^H,{\bf m}^H)$ of $\A$ onto $H$ is defined by 
$$
\A^H:=\{H \cap L \mid L \in \A \setminus \{H\}\},
$$
and by 
$$
{\bf m}^H(X)
:=
|\{L \in \A \setminus \{H\} \mid 
H \cap L=X\}|
$$
for $X \in \A^X$.
\label{Z}
\end{definition}

For a Ziegler restriction $(\A^H,\bf{m}^H)$ in $\K^2 \simeq H$, 
as for $\A$, we may define $$
D(\A^H,\bf{m}^H):=\{
\theta \in \Der \mbox{Sym}^*(H^*) \mid \theta(X) \in 
\mbox{Sym}^*(H^*) \alpha_X\ (\forall X \in \A^H)\}.
$$ Since it is reflexive too and $\dim_\K \K^2=2$, 
it is always free, and we can define the exponents of $(\A^H,\bf{m}^H)$ as for $D(\A)$. Note that when $\bf{m} \equiv 1$, then $D(\A^H,1)=D(\A^H)$, that is free with $\exp(\A^H)=(1,|\A^H|-1)$. 

Note that if $\exp(\A^H,{\bf m }^H)=(d_2,d_3)$ for a Ziegler restriction $(\A^H,\bf{m}^H)$ in 
$\K^2$, then by \cite{Z},
$$
d_2+d_3=|{\bf m}^H|:=
\sum_{X \in \A} {\bf m}^H(X).
$$

\begin{thm}[Theorem 0.3, \cite{A1}]
Let $(\A^H,{\bf m}^H)$ be a balanced 
Ziegler restriction in $\K^2$, i.e., 
there are no $X \in \A^H$ such that 
$$
2{\bf m}^H(X) \ge |{\bf m}^H|:=\sum_{Y \in \A} {\bf m}^H(Y).
$$
Let 
$\exp(\A^H,{\bf m}^H)=(d_1,d_2)$ be the corresponding exponents. Then 
$|d_2-d_1| 
\le |\A^H|-2$.
\label{A1}
\end{thm}

\begin{thm}[\cite{Z}]
Let $\A$ be free with $\exp(\A)=(1,d_2,d_3)$. Then $(\A^H,{\bf m}^H)$ is free with 
$\exp(\A^H,{\bf m}^H)=(d_2,d_3)$.
\label{Zmain}
\end{thm}

The following is maybe well-known, but we include a short proof for the reader's convenience.

\begin{lem}
Let $(\A^H,{\bf m}^H)$ be a Ziegler restriction in $\K^2$ with 
$$|{\bf m}^H|-|\A^H|+1 \le |\A^H|-1.$$ Then $\exp(\A^H,{\bf m}^H)=(|{\bf m}^H|-|\A^H|+1, |\A^H|-1)$.
\label{easy}
\end{lem}

\proof Let 
$$
Q:=\prod_{X \in \A^H} \alpha_X^{{\bf m}^H(X)-1}.
$$
Then for the Euler derivation $\theta_E$ in the derivation module of $\mbox{Sym}^*(H^*)$, it holds that  $Q\theta_E \in D(\A^H,{\bf m}^H)_{|{\bf m}^H|-|\A^H|+1}$. 
Since $D(\A^H,{\bf m}^H)\subset D(\A^H)$ and $\exp(\A^H)=(1,|\A^H|-1)$, 
the larger 
exponent of $(\A^H,{\bf m}^H)$ is at least $|\A^H|-1$. Since 
$|{\bf m}^H|-|\A^H|+1 \le |\A^H|-1$, the smaller exponent of 
$(\A^H,{\bf m}^H)$ is at most $|{\bf m}^H|-|\A^H|+1=\deg Q \theta_E$. 
However, there are no $0 \neq 
\theta \in D(\A^H,{\bf m}^H)$ such that $f\theta=Q\theta_E$ for some $f \in \mbox{Sym}^*(H^*)$. Thus 
$Q\theta_E$ has to be a part of basis, which completes the proof. \endproof

\begin{rk}
In the above we introduced several results for our use. Namely, 
Definition \ref{Z} and Theorem \ref{Zmain} are the three dimensional case of those in \cite{Z}. In \cite{Z} 
these are introduced and proved for 
an arbitrary dimensional arrangement. Also, Theorem \ref{A1} and Lemma \ref{easy} holds true for 
all (multi)arrangements in $\K^2$.
\end{rk}

Now we prove 
Theorem \ref{thm1}, identifying a projective line arrangement in $\PP^2=\PP(\K^3)$ with the corresponding central plane arrangement in $\K^3$.
Let $H \in \A$ and $(\A^H,{\bf m}^H)$ the Ziegler restriction of $\A$ 
onto $H$. First assume that $(\A^H,{\bf m}^H)$ is not balanced for some $H \in \A_p$ with 
a modular point $p$ of maximal multiplicity $m$. Namely, 
there is a modular point $p \in L_2(\A^H)$ such that 
$m=\mult_p(\A) \ge \mult_q(\A)$ for all $q \in L_2(\A)$ by \cite{T1}, and 
${\bf m}^H(p)=m-1 \ge 
\frac{d-1}{2}=\frac{|m^H|}{2}$. In this case $2m-1 \ge d$. Since 
$3m-3 < 2m-1 \iff 
m <2$, in this 
case $d \le 3m-3$.

Thus we may assume that $(\A^H,{\bf m}^H)$ is balanced for all $H \in \A_p$ with a modular point $p$ of 
maximal multiplicity $m$. 
Since $\A$ is free with exponents 
$(1,m-1,d-m)$, it follows that $\exp(\A^H,{\bf m}^H)=(m-1,d-m)$ by Theorem \ref{Zmain}. Assume that 
$d\ge 3m-2$. Then $d-m-(m-1)=d-2m+1=(d-3m+2)+m-1 >0$. Also, since $(\A^H,{\bf m}^H)$ is 
balanced, Theorem \ref{A1} shows that 
$$
d-2m+1 \le |\A^H|-2 \iff
|\A^H| \ge d-2m+3
\ge m+1.
$$
So for $\exp(\A^H)=(1,|\A^H|-1)$, it holds that 
$|\A^H|-1 \ge m$. Since for any $\mu:\A^H \rightarrow \Z$ with $1 \le \mu \le m^H$, by Lemma 
\ref{easy}, it holds that 
$\exp(\A^H,\mu)=(|\mu|-|\A^H|+1,|\A^H|-1)$ if $|\mu|-|\A^H|+1 \le |\A^H|-1$. Since $m-1 
<|\A^H|-1$, $\exp(\A^H,m^H)=(m-1,d-m)$ occurs only when $d-m=|\A^H|-1 
\iff |\A^H|=d-m+1$. Since $p$ is a point on $H$ 
with multiplicity $m$, it holds that $|\A_q|=2$ for all $q \in L_2(\A)$ with $p \neq q \in H$. Namely, for 
all $H \in \A_p$ and $L \in \A \setminus \A_p$, $|\A_{H \cap L}|=2$, but this cannot occur for all such 
$H$ unless $d-m=|\A \setminus 
\A_p|=1$ since $\A$ is supersolvable and $p$ is a modular point. So $d \le 3m-3$. 
The  proof of Theorem \ref{thm1} is completed.

\section{Proof of Theorem \ref{thm1B}} 

Assume that  two of the modular points are $p=(0:0:1)$ and $p'=(0:1:0)$, and hence the line
$L_0: x=0$ is in the arrangement $\A$. Assume that the $(m-1)$ lines in $\A$ passing through $p$, and different from $L_0$, are given by the equations
$$L_j: y-b_jx=0 \text{ for } j=1,2, \ldots m-1,$$
for some distinct complex numbers $b_j$. Similarly, the $(m-1)$ lines in $\A$ passing through $p'$, and different from $L_0$, are given by the equations
$$L'_j: z-c_jx=0 \text{ for } j=1,2, \ldots m-1,$$
for some distinct complex numbers $c_j$. Note that one has
$$L_i\cap L'_j=(1:b_i:c_j)$$
for any pair $(i,j)$. 

If $d=|\A|=2m-1$, then we have listed all the lines in $\A$, and the proof is complete, namely $\tilde \A$ is lattice-isotopic to $\tilde \A(w_0)$, since the Zariski open set
of all vectors $a=(a_1,\ldots,a_{m-1})\in \C^{m-1}$ with $a_i \ne a_j$ for $i \ne j$ is path-connected by smooth arcs. 

If $d>2m-1$, then let $L$ be a line in $\A$, different from the $2m-1$ lines already listed above. Then there is an index $i$ such that $q_1=L \cap L_1=L \cap L'_i$.
We can and do assume that $i=1$ and $b_1=c_1=0$, in other words $q_1=(1:0:0)$.
Then, for any $j=2,\ldots,m-1$, there is an index $\sigma(j) \in \{2,\ldots,m-1\}$ such that
$$q_j=L\cap \L_j=L \cap L'_{\sigma(j)}=(1:b_j:c_{\sigma(j)}).$$
In this way we get a permutation $\sigma$ of the set $\{2,\ldots,m-1\}$.
The points $q_1,\ldots,q_{m-1}$ are collinear if and only if the two $(m-2)$ vectors $u=(b_2, \ldots, b_{m-1}) \ne 0$ and $v=(c_{\sigma(2)}, \ldots, c_{\sigma(m-1)}) \ne 0$ are proportional, i.e. there is $\lambda \in \C^*$ such that $v=\lambda u$.
An equation for the union of lines $L_1 \cup \ldots \cup L_{m-1}$ is
$$P(x,y)=\prod_{j=1,m-1}(y-b_jx),$$
and an equation for the union of lines $L'_1 \cup \ldots \cup L'_{m-1}$ is
$$P'(x,z)=\prod_{j=1,m-1}(z-c_jx).$$
Under the coordinate change $(x,y,z) \mapsto (x,y,Z)$, where
 $z=\lambda Z$, the polynomial $P$ becomes
$P(x,y)=yP_0(x,y)$ and the polynomial $P'$ becomes
$P'(x,Z)=\lambda^{m-1}ZP_0(x,Z)$, where
$$P_0(s,t)=\prod_{j=2,m-1}(t-b_js),$$
is a polynomial with $(m-2)$ distinct linear factors.
Under this change of coordinates the point $q_1$ remains $(1:0:0)$.
To simplify the notation, we denote the new coordinates by $(x,y,z)$. The conclusion is that, if $d>2m-1$, we may assume that the constants $b_j $ and $c_j$ above satisfy $b_1=c_1=0$ and $b_j=c_j$ for $j>1$.

The line $L$ is determined by the two points $q_1=(1:0:0)$ and $q_2=
(1:b_2:b_{\sigma(2)})$, in other words it is determined by $\sigma(2)$, and has the equation
$$L:b_2z-b_{\sigma(2)}y=0.$$
As explained above, one has
$$v=(b_{\sigma(2)}, \ldots, b_{\sigma(m-1)}) =\lambda (b_2, \ldots, b_{m-1})=\lambda u.$$
Let $n=m-2$ as in Introduction.
By taking the product of the $n$ components of the vectors $u$ and $v$ we get 
$$\prod_{j=2,m-1}b_j=\lambda^{n}\prod_{j=2,m-1}b_j,$$
and hence 
\begin{equation}
\label{LA}
\lambda^{n}=1, 
\end{equation}
since $\prod_{j=2,m-1}b_j \ne 0$. It follows that, for the given line $L$, there is a root of unity $\lambda _L\in \mu_n$ such that 
$$L:z-\lambda_Ly=0.$$
Let $e_L$ be the order of the root of unity $\lambda _L$, and let $s_j$ be the $j$-th symmetric function in $b_2, \ldots, b_{m-1}$. Then as above we get
$$s_j=\lambda^{j}s_j,$$
and hence $s_j=0$ when $j$ is not a multiple of $e_L$.
This implies that the polynomial $P_0(s,t)$ above has the special form
$$P_0(s,t)=Q_0(s^{e_L},t^{e_L}),$$
for some homogeneous polynomial $Q_0(s,t)$ of degree $n'=n/e_L$, having $n'$ distinct linear factors.
In conclusion, the number $k$ of lines $L$ in $\A$, through the point $q_1$, and distinct from $L_1$ and $L'_1$, satisfies 
$$0 \leq k \leq n=m-2.$$

\medskip

\noindent {\bf Case 1: all the lines $L$ in $\A$, not passing through one of the two modular points $p$ and $p'$, have a common intersection point $q$.}

This case is covered by the discussion above, since we can assume that $q_1$ is the common intersection point $q$.  We can assume that $\A$ is  given by the equation
$$\A: xyzQ_0(x^e,y^e)Q_0(x^e,z^e)\prod_{j=1,k}(z-\lambda_{L^j} y)=0,$$
with $\lambda_{L^j}^e=1$ and the product being taken over the $k$ lines 
$L^1$, \ldots, $L^k$ in $\A$ passing through $q_1$ and different from $L_1$ and $L_2$. Here $e$ is the least common multiple of the orders of the various $\lambda_{L^j}$ in the above product. If we set 
$$w=(\lambda_{L^1},\ldots, \lambda_{L^k}),$$
to show that the central plane arrangements $\tilde A$ and
$\widetilde {A(w)}$ are lattice-isotopic, it is enough to deform the polynomial
$Q_0(s,t)$ to the polynomial $Q'_0(s,t)=s^{n'}-t^{n'}$ within the class of 
degree $n'$ binary forms in $s,t$ without multiple factors. This set of binary forms, being the complement of the corresponding discriminant hypersurface, is clearly path-connected by smooth arcs. Note that when $k=n$, then $e=n$ and hence $Q_0(s,t)=as+bt$ for some $a,b \in \C^*$, being a linear form. This fact implies that in this case, $\A$ is projectively equivalent to  $\A(m-2,1,3)$.
This completes the proof in this case, except for the uniqueness of $w$ which is proved in the next Lemma.

\medskip

\noindent {\bf Case 2: all the lines $L$ in $\A$, not passing through one of the two modular points $p$ and $p'$, do not have a common intersection point.}

We show that this case is not possible. By the above discussion, we can assume that $\A$ has a subarrangement $\B$, given by the equation
$$\B: xyzQ_0(x^e,y^e)Q_0(x^e,z^e)\prod_{j=1,k}(z-\lambda_{L^j} y)=0,$$

Since we are in Case 2,  there is at least one new line $L' \in \A \setminus \B$, with an equation
$$L': x+by+cz=0,$$
where $bc \ne 0$. This line creates the following new intersection points.

\medskip

\noindent (i) The point $p''=L \cap L'$, where we choose one line $L=L^1$ as in the discussion above. Then $p''=(-b-c\lambda_L:1:\lambda_L)$.
This point can be either on $L_0$, and then $b+c \lambda_L=0$, or
if $b+c \lambda_L \ne 0$, there is an index $j_0$ and an index $k_0$ such that $p''=L_{j_0} \cap L'_{k_0}$. This latter equality holds if and only if
$$1+b_{j_0}(b+c \lambda_L)=0 \text{ and } \lambda_L+b_{k_0}(b+c \lambda_L)=0.$$

\medskip

\noindent (ii) The points $p_j=L_j \cap L'=(c:cb_j:-(1+bb_j))$ for $j=2,\ldots,m-1$. As above, there is a permutation $\sigma'$ of the set $\{2,\ldots,m-1\}$ such that $p_j \in L'_{\sigma'(j)}$. This means
$$1+bb_j+cb_{\sigma'(j)}=0.$$
If we sum up these equalities for $j=2,\ldots,m-1$ and denote $K=\sum_{j=2,\ldots,m-1}b_j$, we get
$$m-2+(b+c)K=0.$$
Note that when $e>1$, we get $K=0$ and hence a contradiction, since $m \geq 3$. It follows that the only possible case is when $e=1$, hence there is a unique line $L:y-z=0$ through $q_1$, and passing through none of the two modular points $p$ and $p'$.
Hence, the only possible case is  $e=1$ and $K \ne 0$, which implies that the sum
$$b+c=(2-m)K^{-1}=:K'$$
is uniquely determined. The equation of $L'$ can then be written as
$$L'=L'_b: x+by+(K'-b)z=x+K'z +b(y-z)=0,$$
and hence depends only on the parameter $b$.
Note that the point $p''$ is now given by $(-K':1:1)$, and hence, when $b$ varies, all the lines  $L_b' \in \A \setminus \B$ pass through $p''$.
Since $L$ also passes through $p''$, we have reached a contradiction,
namely the arrangement $\A$ is in  Case 1 above and not in Case 2.

\begin{lem}
\label{lemU}
The isomorphism class of the lattice $L(w)$ determines the class 
$[w] \in \overline{ W(n,k)}$ uniquely.
\end{lem}
\proof
The cases $k=0$, $k=1$ and $k=n$ are obvious, since $\overline{ W(n,k)}$ are singletons in these cases. We assume that $2 \leq k\leq n-1$ and hence $m>k+2=4$. We use the notations in the above proof. The arrangement $\A(w)$ has two modular points of multiplicity $m$, namely $p$ and $p'$,
 a point $q=q_1$ of multiplicity $k+2$ with  $4 \leq k+2<m$, and other intersection points of multiplicity 2 and 3. We can and do assume that
 $$w=(1,w_2,\ldots,w_k).$$
The fact that the line $L^1:y-z=0$ is present, allows us to parametrize the $n$ lines
$L_i$ through $p$ avoiding $p'$ and $q$, and the $n$ lines $L'_j$ through $p'$ avoiding $p$ and $q$, using as index set the group $\mu_n$,
in such a way that $L^1 \cap L_a \cap L'_{a'} \ne \emptyset$ for $a,a' \in \mu_n$ if and only if
$a=a'$. Any other line $L^j$ for $j=2,\ldots, k$ through $q$ avoiding $p$ and $p'$
gives rise to a bijection $\sigma^j:\mu_n \to \mu_n$, which was shown above to be the translation by an element $w_j \in \mu_n$. In this way we get the 
unordered set $\{w_2,...,w_k\}$, which yields a well defined class
$$[w]=[(1,w_2,\ldots,w_k)] \in \overline{ W(n,k)}.$$
This construction is canonical, up-to the interchange of the two modular points $p$ and $p'$. This operation changes $w$ into
$$w'=(1,w_2^{-1},\ldots,w_k^{-1})$$
and clearly $[w']=[w]$. 
\endproof

\begin{rk}
\label{rkHH}
If the line arrangement $\A$ in Theorem \ref{thm1B} is real, then in the above proof we see that equation \eqref{LA} implies that either $\lambda=1$, or $\lambda=-1$ and $m\geq 4$ even. So for real line arrangements, only the cases $k=0,1$  and any $m$, and $k=2$ and $m\geq 4$ even can occur, exactly as explained already in \cite[Section 3.2.3]{HaHa}. To see that the case $k=2$ and $m\geq 4$ even is realizable from our point of view, just take the set of $m-2$ nonzero real numbers
$\{b_2,\ldots, b_{m-1}\}$ in the above proof to be invariant under multiplication by $-1$.
\end{rk}

\begin{rk}
\label{rkMV}
Various other subarrangements of the  full monomial line arrangement 
$\A(n,1,3)$ have been studied by Marchesi and Vall\`es in \cite{MV},
in relation with free line arrangements, not necessarily supersolvable,  and Terao's conjecture.
\end{rk}

\section{Proof of Theorem \ref{thm2B}} 

 The numbers $n_k$, of the intersection points of multiplicity $k$ in a line arrangement $\A$, satisfy a number of relations. The easiest of them is  the following.
	\begin{equation}
\label{eqSum}
\sum_{k\geq 2} n_k{k \choose 2}={d \choose 2}, 
\end{equation}
where $d=|\A|$.
A highly non-trivial restriction on these numbers is given by the Hirzebruch inequality, valid for non trivial line arrangements, see \cite{H}:
\begin{equation}
\label{eqHir}
n_2 +\frac{3}{4}n_3 -d\geq  \sum_{k >4}(k-4)n_k.
\end{equation}

Moreover, if $\A$ is supersolvable, then it is shown in \cite[Proposition 3.1]{T1} that the following holds
\begin{equation}
\label{eqSS}
n_2  \geq 2\sum_{k\geq 2}n_k-m(d-m)-2,
\end{equation}
where $m=\max\{k : n_k \ne 0\}$. 
\begin{rk}
\label{rkEQ}
It is shown in \cite[Proposition 3.1]{T1} that the equality holds in \eqref{eqSS} if and only if the supersolvable line arrangement $\A$ 
is a cone $C(\A')_e$ over a generic line arrangement $\A'$, as defined in Introduction. See Theorem \ref{thm2C} for more on  such cones.
\end{rk}

Denote by $\tau(\A)$ the global Tjurina number of the arrangement $\A$, which is the sum of the Tjurina numbers $\tau(\A,a)$ of the singular points $a$ of $\A$. One has
\begin{equation} \label{r0}
\tau(\A)= \sum_{k \geq 2}n_k(k-1)^2.
\end{equation}
Indeed, any singular point $a$ of multiplicity $k \geq 2$ of a line arrangement $\A$ being weighted homogeneous, the local Tjurina number $\tau(\A,a)$ coincides with the local Milnor number $\mu(\A,a)=(k-1)^2$, see \cite{KS}.

At this point we recall some notation.
Let $S=\C[x,y,z]$ be the polynomial ring in three variables $x,y,z$ with complex coefficients, and let $\A:f=0$ be an arrangement of  $d$ lines in the complex projective plane $\PP^2$. The minimal degree of a Jacobian syzygy for the polynomial $f$ is the integer $mdr(f)$
defined to be the smallest integer $r \geq 0$ such that there is a nontrivial relation
\begin{equation}
\label{rel_m}
 af_x+bf_y+cf_z=0
\end{equation}
among the partial derivatives $f_x, f_y$ and $f_z$ of $f$ with coefficients $a,b,c$ in $S_r$, the vector space of  homogeneous polynomials in $S$ of degree $r$. This is the same as the minimal degree of a non-zero derivation $\theta \in D(\A)$ such that $\theta f=0$.
It is known that  $mdr(f)=0$ if and only if $n_d=1$, hence  $\A$ is a pencil of $d$ lines passing through one point. Moreover, 
$ mdr(f)= 1$ if and only if $n_d=0$ and $n_{d-1}=1$, hence  $\A$ is a 
near pencil in the terminology of \cite{HaHa}. A line arrangement $\A:f=0$ is called nontrivial, following \cite{HaHa}, if $mdr(f) \geq 2$.
Note that the line arrangements $\A:f=0$ satisfying $mdr(f)=2$ are classified in \cite{To}. As said in Section 2, any nontrivial supersolvable line arrangement $\A:f=0$ is free, with exponents $\exp(\A)=(1,d_2,d_3)$, where
$d_2=r$ and $d_3=d-1-r$, with $r=mdr(f)$ and $d=|\A|$.
In particular, using \eqref{FA2}, we see that for a free, and hence for a
supersolvable arrangement, one has $mdr(f) \leq (d-1)/2$.

\begin{rk}
\label{rkm=3}

Let $\A$ be  a $3$--homogeneous supersolvable line arrangement, which is not a pencil.
We get $d=|\A] \leq 6$ from Theorem \ref{thm1}. 
The 3-homogeneous supersolvable line arrangements with $n_3>1$ are classified in \cite[Section (3.2.2)]{HaHa}. Here we show that the case
$n_3=1$ is not possible, and give an alternative approach for the case $n_3>1$.
 By the above discussion, we have two possibilities.
\begin{enumerate}
\item $d=6$ and then $mdr(f)=m-1=2$. Using the classification in \cite{To} of line arrangements with $mdr(f)=2$, 
we see that the arrangement $\A$ is linearly equivalent to the full monomial line arrangement $\A(1,1,3)$ above. For this arrangement
$n_2=3=d/2$.

\item $d=5$ and then again $mdr(f)=2$. By the classification in \cite{To}, 
we see that the arrangement  $\A$ is linearly equivalent to the arrangement
$$ \B: xyz(x-y)(x-z)=0,$$
obtained from $\A(1,1,3)$ by deleting one line, namely $y-z=0$. For this arrangement
$n_2=4>d/2$. 

\end{enumerate}
\end{rk}

  Note that, for any  free line arrangement with exponents $\exp(\A)=(1,d_2,d_3)$, one has
$$\tau(\A)=(d-1)^2-d_2d_3,$$
see for instance \cite{DSer}.
In particular, for an $m$--homogeneous supersolvable line arrangement, consisting of $d=|\A|$ lines one has
\begin{equation} \label{tauSS}
\tau(\A) = (d-1)^2-(m-1)(d-m).
\end{equation}
In fact, {\it this formula holds more generally} for any free line arrangement $\A$, if we assume $2m \geq d$, where $m$ is the maximal multiplicity of an intersection point in $\A$. To see this, we use \cite[Theorem 1.2]{Dcurves}, which shows that the condition $2m \geq d$ implies that either $d_2=mdr(f)=d-m, d_3=m-1$, or $d_2=mdr(f)=m-1$ and $d_3=d-m$.
 The formulas \eqref{r0} and \eqref{tauSS} imply
\begin{equation}
\label{eqSS1}
\sum_{k\geq 2}n_k(k-1)^2=(d-1)^2-(m-1)(d-m).
\end{equation} 

\medskip

\noindent { \bf Proof of Theorem \ref{thm2B}}

\medskip

We use first the relations \eqref{eqSum} and \eqref{eqSS1}. By eliminating $n_3$ we get
$$n_2=2d(d-1)-3\tau(\A)+\sum_{k \geq 4}n_k(k^2-4k+3) \geq 2d(d-1)-3\tau(\A)+m^2-4m+3.$$
If we use the formula \eqref{eqSS1}, we get
$$n_2 \geq -2m^2+(3d-1)m-d^2+d.$$
Hence it is enough to show that $2m \geq d$ implies
$$-2m^2+(3d-1)m-d^2+d \geq \frac{d}{2}.$$
This is equivalent to
$$-4m^2+2(3d-1)m-2d^2+d=-(2m-d)(2m-2d+1) \geq 0,$$
which is clear since $2m-d \geq 0$ by our hypothesis, and $2m-2d+1<0$ since $\A$ is not a pencil.
This ends the proof of Theorem \ref{thm2B}.

\section{Proof of Theorem \ref{thm2C}, Proposition \ref{prop2C} and Corollary \ref{corthm2}.} 

Now we prove Theorem \ref{thm2C}, using the notation in Introduction, when we have defined the cone $C(\A')_e$.
Note that passing from $C(\A')_0$ to $C(\A')_e$ increases
the multiplicity of the vertex $p$ and the total number of lines by $e$,
while the number of double points is increased by $ed'$. Hence it is enough to consider the case $e=0$.
First we discuss the case $d'=|\A'|=2d_1+1$ is odd.
 Then the total number of double points in $\A'$ is $N=d_1(2d_1+1)$. By Bezout Theorem, any line in $\PP^2$, not in $\A'$, can contain at most $d_1$ among these $N$ double points. It follows that $N' \geq d'$, and this is equivalent to
Conjecture \ref{conj1} for $\A=C(\A')$, since clearly $d=d'+N'$ and $$n_2=d'(N'-d'+1),$$ by \cite[Proposition 3.1]{T1} quoted in Remark \ref{rkEQ}.

When $d'=2d_1$ is even,  then $N=d_1(2d_1-1).$
 Conjecture \ref{conj1} holds for $\A$ with the same argument as above, unless we are in the following {\it exceptional case:
 $m=N'=2d_1-1$ and any line determined by $p$ and one of the $N$ nodes contains exactly $d_1$ nodes. 
 We show that this case cannot occur. This is a special case of \cite[Conjecture 3.3]{T1}, since any line coming from $\A'$ would contain only triple points in $\A$ in this exceptional case.} 
 
 Let $g=0$ be an equation for the generic arrangement $\A'$, and let $J_g$ be the Jacobian ideal of 
$g$, spanned by $g_x,g_y,g_z$ in $S=\C[x,y,z]$. Then the minimal resolution of the Milnor algebra $S/J_g$ has the following form
$$0 \to S^{d'-3}(2-2d') \to S^{d'-1}(3-2d') \to S^3(1-d')\to S,$$
see \cite{RT,Yu}. Using \cite[Proposition 1.3 (i)]{HS}, we see that the quotient $J_g^{sat}/J_g$, 
where $J_g^{sat}$ is the saturated ideal of $J_g$ with respect to the maximal ideal $(x,y,z)$,
has $d'-3$ generators as a graded $S$-module, all in degree $d'-1$.
It follows that $\dim (J_g^{sat})_{d'-1}=d'$. Since $\A'$ has only nodes as singularities, it follows that $J_g^{sat}=\sqrt J_g$, and hence a homogeneous polynomial $h \in J_g^{sat}$ if and only if $h$ vanishes on the set of $N$ nodes.
Let $L_1,...,L_{d'}$ be the lines in $\A'$ and choose $\ell_i=0$ an equation for the line $L_i$, for $i=1, \ldots, d'$.
Note that 
$$g=\prod_{i=1,d'} \ell_i,$$
and the following set of $d'$ homogeneous polynomials of degree $d'-1$
$$g_i= \frac{g}{\ell_i}$$
for $i=1, \ldots,d'$ is contained in $(J_g^{sat})_{d'-1}$. Since these polynomials are clearly linearly independent, they form a basis of the $\C$-vector space $(J_g^{sat})_{d'-1}$.
{\it Assume now that we are in the exceptional case above.}
Choose the coordinates such that the vertex $p=(0:0:1)$, and assume that the line $y=0$ {\it is not}  among the $m=d'-1$ lines passing through $p$. Note that in these coordinates, we can choose
$$\ell_i=z-h_i(x,y) \text{ with } h_i(x,y)=a_ix+b_iy$$
where $a_i,b_i \in \C$, and the $a_i$'s are all {\it distinct}.
In fact, the union of the $m$ lines through $p$ are given by
any of the following equations
$$H_i=\prod_{j=1,\ldots, d', j\ne i}(h_i(x,y)-h_j(x,y))=\prod_{j=1,\ldots, d', j\ne i}((a_i-a_j)x+(b_i-b_j)y)=0,$$
for $ i=1, \ldots,d'$, which differ from each other by nonzero constant factors. 
Since $H_1$ vanishes at all the $N$ nodes of $\A'$, it follows that
$H_1 \in (J_g^{sat})_{d'-1}$, and hence there are constants $\lambda_i \in \C$ such that
\begin{equation}
\label{key}
\sum_{i=1,d'}\lambda_ig_i=H_1.
\end{equation}
We can write 
$$g_i=\prod_{j=1,\ldots, d',j\ne i}(z-h_j)=z^{d'-1}-\sigma_1^iz^{d'-2}+\sigma_2^iz^{d'-3}+\ldots +(-1)^{d'-1}\sigma_{d'-1}^i,$$
where $\sigma_k^i$ denotes the $k$-th symmetric function of $h_1, \ldots, \widehat h_i, \ldots, h_{d'}$, with $\sigma_0^i=1$ by convention. If we denote by $\sigma_k$ the $k$-th symmetric function of $h_1, \ldots,  h_i, \ldots , h_{d'}$, with $\sigma_0=1$, then one has the following relations
\begin{equation}
\label{sym}
\sigma_k^i=\sigma_k-\sigma_{k-1}h_i+\sigma_{k-2}h_i^2- \ldots+(-1)^k\sigma_0h_i^k.
\end{equation}
Using the equation \eqref{key} we get the following system of equations:
$$(E_k): \  \  \  \sum_{k=1,d'}\lambda_i\sigma_k^i=0,$$
for $k=0,\ldots,d'-2$.
Using now the equation \eqref{sym}, the system $(E_i)$ can be transformed into the following equivalent system of equations:
$$(E'_k): \  \  \  \sum_{i=1,d'}\lambda_ih_i^k=0,$$
for $k=0,\ldots,d'-2$.
If we consider the coefficient of $x^uy^v$, for $u+v=k$, in the equation $(E'_k)$, we get the following system of equations
$$(E''_{u,v}): \  \  \  \sum_{i=1,d'}\lambda_ia_i^ub_i^v=0.$$
for $0 \leq u+v \leq d'-2$, with $u,v \geq 0$. Since $a_i \ne a_j$ for $i \ne j$, it follows that the vectors 
$${\bf a}^{(j)}=(a_1^j,a_2^j, \ldots, a_{d'}^j),$$
for $j=0,1,\ldots, d'-1$ are linearly independent. Consider also the vector
$${\bf b}=(b_1,b_2, \ldots, b_{d'}),$$
and note that the vectors ${\bf a}^{(0)}, {\bf a}^{(1)}$ and ${\bf b}$ are linearly independent since the line arrangement $\A'$ has only double points. 

\medskip

\noindent {\bf Case 1}. If the vectors ${\bf a}^{(j)}=(a_1^j,a_2^j, \ldots, a_{d'}^j),$
for $j=0,1,\ldots, d'-2$ and $\bf b$ are linearly independent, then the above equations $(E''_{u,v})$ in $\lambda_i$ for the pairs $(u,v)=(0,0), (1,0), \ldots, (d'-2,0)$ and respectively $ (0,1)$, imply that $\lambda_i=0$ for all $i$, a contradiction.

\medskip

\noindent {\bf Case 2}. If the vectors ${\bf a}^{(j)}=(a_1^j,a_2^j, \ldots, a_{d'}^j),$
for $j=0,1,\ldots, d'-2$ and $\bf b$ are linearly dependent, then there are unique coefficients $\al_i \in \C$ such that
$${\bf b}=\sum_{j=0,1,\ldots, d'-2}\al_j{\bf a}^{(j)}.$$
Let $M=\max \{j : \al_j \ne 0\}$ and note that $2 \leq M \leq d'- 2$.
If we look coordinate-wise, we have
$$b_i=\sum_{j=0,1,\ldots, d'-2}\al_ja_i^j.$$
We multiply this equation by $a_i^{k-1}$, with $2\leq k=d'-M \leq d'-2$
and get
$$a_i^{k-1}b_i=\sum_{j=0,1,\ldots, M}\al_ja_i^{j+k-1},$$
where the last term in the sum is $\al_Ma_i^{d'-1}$.
It follows that the vector
$${\bf c}=(a_1^{k-1}b_1, a_2^{k-1}b_2, \ldots, a_{d'}^{k-1}b_{d'})=\sum_{j=0,1,\ldots, M}\al_j{\bf a}^{(j+k-1)}$$
is linearly independent of the vectors ${\bf a}^{(j)}=(a_1^j,a_2^j, \ldots, a_{d'}^j),$
for $j=0,1,\ldots, d'-2$.  Hence the above equations $(E''_{u,v})$ for the pairs $(u,v)=(0,0), (1,0), \ldots, (d'-2,0)$ and respectively $ (k-1,1)$, imply that $\lambda_i=0$ for all $i$, a contradiction with the equation \eqref{key}. This completes the proof of Theorem 
\ref{thm2C}.

To prove Proposition \ref{prop2C}, we note as above that it is enough to treat the case $e=0$. Let $n'_k$ be the number of points of multiplicity $k$ in $\A'$.
If $\A'$ is a near-pencil, that is $n'_{d'-1}\ne 0$, then the total number of singular points in $ \A'$ is $N'=d'$. Hence $d=|\A|=2d'$ and $m=m_p(\A)=d'$, which shows that Conjecture \ref{conj1} holds in this case, since $d \leq 2m$ and we apply Theorem \ref{thm2B}.
Assume from now on that
$\A'$ is neither a pencil, nor a near-pencil. Then Hirzebruch's inequality
\eqref{eqHir} holds, and it implies
$$N' \geq n_2'+n_3' \geq d'.$$
We have as above $d=d'+N'$ and $m=m_p(\A)=N'$, which shows that Conjecture \ref{conj1}  holds, since $d \leq 2N'$ and we apply again Theorem \ref{thm2B}.

Finally we prove Corollary \ref{corthm2}. 
Note that Theorem \ref{thm1} implies that for $m\leq 5$ we get $d \leq 12$. The claim then follows from \cite[Corollary 3.5]{T1}.

\end{document}